\documentclass[11pt]{article}
\def\today{18/05/08} 
\usepackage{amsmath,amsfonts,amssymb}
\usepackage{amsthm,amscd}

\newtheorem{theorem}{Theorem}[section]
\newtheorem{lemma}{Lemma}[section]

\newtheorem{corollary}[theorem]{Corollary} \theoremstyle{definition}
 \theoremstyle{remark}
\newtheorem{remark}{Remark}[section]

\newcommand{\R}{{\mathbb R}} \newcommand{\U}{{\mathcal U}}
\newcommand{\V}{{\mathcal V}} \newcommand{\A}{{\mathcal A}}
\newcommand{\F}{{\mathcal F}} \newcommand{\G}{{\mathcal G}}
\newcommand{\Lc}{{\mathcal L}} \newcommand{\Hc}{{\mathcal H}}
\newcommand{\Z}{{\mathbb Z}} 
 \newcommand{\Jc}{{\mathcal J}}
\newcommand{\Ic}{{\mathcal I}} \newcommand{\Oc}{{\mathcal O}}
\newcommand{\Sc}{{\mathcal S}} \newcommand{\T}{\mathbb{T}}
\newcommand{\bn}{{\bf n}}

\def\Nc{{\mathcal N}}
\def\M{{\mathcal M}}
\def\be{{\bf e}_1}

\def\norma#1{\left \| #1\right\|}
\def\Ph{{\mathcal P}}

\numberwithin{equation}{section}

\begin{document}

\author{D. Bambusi, C. Bardelle \\ {bambusi,bardelle@mat.unimi.it}, \\
Universit\`a degli Studi di Milano, \\ Dipartimento di Matematica
``F. Enriques'', \\ Via Saldini 50, 20133 Milano, Italy}

\title{Invariant tori for commuting Hamiltonian PDEs}

\date{\today}

\maketitle

\begin{abstract}
We generalize to some PDEs a theorem by Nekhoroshev on the
persistence of invariant tori in Hamiltonian systems with $r$
integrals of motion and $n$ degrees of freedom, $r\leq n$. The
result we get ensures the persistence of an $r$-parameter family of
$r$-dimensional invariant tori. The parameters belong to a
Cantor-like set. The proof is based on the Lyapunof-Schmidt
decomposition and on the standard implicit function theorem. Some of
the persistent tori are resonant.  We also give an application to
the nonlinear wave equation with periodic boundary conditions on a
segment and to a system of coupled beam equations. In the first case
we construct 2 dimensional tori, while in the second case we
construct 3 dimensional tori.
\end{abstract}


\setcounter{section}{0}

\section{Introduction}
\label{S.0}

In 1994 Nekhoroshev \cite{N94} published a theorem on persistence of
invariant tori in Hamiltonian systems with symmetry. Such a theorem
interpolates between Poincar\'e's theorem on persistence of periodic
orbits and Arnold Liouville theorem on the structure of the phase
space of integrable systems.

More precisely Nekhoroshev's theorem pertains a system of $r$
commuting Hamiltonians $H^{(1)},...,H^{(r)}$ in a $2n$ dimensional
phase space, $n\geq r$; assume (i) that there exists one
$r$-dimensional torus $\T^r_0$ invariant under the flow of each of the
Hamiltonians $H^{(l)}$ and (ii) that a suitable nonresonance condition
is fulfilled. Then Nekhoroshev's theorem ensures that there exists a
$2r$-dimensional symplectic manifold $N$ on which each system is
integrable; in particular $N$ is foliated in $r$-dimensional tori
which are invariant under the flow of each system. In \cite{BG02} a complete proof of Nekhoroshev's
theorem was published, and the result was extended to the case where
the Hamiltonians $H^{(l)}$ are slightly perturbed. An application to
some infinite dimensional systems was given in \cite{BV02} where the
existence of quasiperiodic breathers in some chains of particles with
symmetry was proved.

In the present paper we extend Nekhoroshev's theorem to a suitable
class of Hamiltonian PDEs and give applications to a nonlinear Klein
Gordon equation with periodic boundary conditions on a segment and to
a system of coupled beam equations. In the first case we construct 2
dimensional tori, while in the second case we construct 3 dimensional
tori.

The extension is far from trivial since the case of PDEs involves
small denominators. To overcome such a problem we have to impose (i) a
quantitative nonresonance condition that generalizes the assumption of
Nekhoroshev's theorem and (ii) a further non degeneracy condition
similar to that used in KAM theory. Moreover, the proof we give here
is completely different from the one of \cite{BG02} which was based on
the use of generalized Poincar\'e sections. On the contrary the
present proof is based on the Lyapunof--Schmidt decomposition and on
the use of a strong nonresonance condition (see \cite{bam00}) which
allows to avoid KAM-type techniques. Indeed we only use the standard
implicit function theorem.

{Related results were obtained in the papers
\cite{Pro05,Bal06,BerPro06} in which the authors exploites translation
invariance in order to construct some quasiperiodic solutions of the
nonlinear wave equation. While the theory of
\cite{Pro05,Bal06,BerPro06} is based on the structure of the wave
equation, the theory we develop here applies to the general
case of cummuting Hamiltonian systems, including the case where all
the Hamiltonians are nonlinear. Our theory however does not apply
directly to the wave equation treated in \cite{Pro05,Bal06,BerPro06}
since in that case our nonresonance assumption is violated.}

We also recall that persistence of invariant tori in systems with
symmetries is often studied by passing to the reduced system and
then using continuation techniques for periodic orbits of such a
reduced system (see e.g. \cite{Mon97,COR02}). Such techniques do not
seem to be applicable to the case of PDEs where often the orbits of
the symmetry group are continuos but not differentiable (they are
orbits of a partial differential equation).

Persistence of invariant tori in Hamiltonian PDEs is usually
obtained through KAM techniques (see e.g.
\cite{K87,CW92,CY00,Bo98,Bou03,EK06}). The main difference with the
previous papers is that we deal here with systems with symmetry.
This allows to obtain persistence of invariant tori using just the
standard implicit function theorem based on the contraction mapping
principle. Moreover the flow on the tori whose existence is proved
here is not necessarily ergodic. Actually, in the case of the wave
equation we prove the existence of some tori on which the flow is
rational.

The paper is organized as follows. In section \ref{stat} we state the
main theorem of the paper, namely the generalization to some PDEs of
Nekhoroshev's theorem. In sect. \ref{proof} we give the proof of the
main theorem. Such a section is split into a few subsections. Finally
in section \ref{Appli} we give the application to the nonlinear wave
and the nonlinear beam equation.

\noindent{Acknowledgements} DB would like to thank Tom Bridges, for a
very stimulating discussion on the families of invariant tori in the
nonlinear wave equation. Such a discussion was the starting point of
the present work.

\section{Statement of the main result}\label{stat}

Fix $r\geq 1$; consider the spaces $\ell^2_{s,\sigma}$ of the
sequences $p=\{ p_j\}_{j\geq r+1}$ s.t.
\begin{eqnarray}
\label{e.1.4}
\norma{p}^2_{s,\sigma}:=\sum_{j}[j]^{2s}{\rm
e}^{2\sigma j}\left|\hat p_j \right|^2<\infty\ ,\quad
[j]:=\max\left\{1,j\right\}\ ,
\end{eqnarray}
and define the phase spaces as
\begin{equation}
\label{e.2}
\Ph_{s,\sigma}:=\U\times \T^r\times \ell^2_{s,\sigma}\times
\ell^2_{s,\sigma}\ni (I,\phi,p,q)\equiv z\
\end{equation}
where $\U\subset\R^r$ is open, $\T=\R/\Z$. We endow $\Ph_{s,\sigma}$
with the norm (on the tangent space of
$\Ph_{s,\sigma}$)
\begin{equation}
\label{e.e.2}
\norma{z}_{s,\sigma}^2=\norma I^2+\norma \phi^2+\norma p_{s,\sigma}^2
+\norma q_{s,\sigma}^2\ .
\end{equation}
We will also use the weak scalar product
\begin{equation}
\label{sca.1}
\langle (I,\phi,p,q);(I',\phi',p',q')
\rangle:=\sum_{j=1}^r\left(I_jI'_j+\phi_j\phi'_j\right) +\sum_{j\geq
  r+1} \left(p_jp_j'+q_jq'_j\right)\ .
\end{equation}
The open ball of radius $R$ and center $0$ in $\R^r$, $\Ph_{s,\sigma}$
or $\ell_{s,\sigma}^2$ will be denoted by $B_R$.  Define the Poisson
tensor $\Jc$ by
\begin{equation}
\label{poi.0}
\Jc(I,\phi,p,q)=(-\phi,I,-q,p )\ .
\end{equation}
Then, given a Hamiltonian function $H$ on $\Ph_{s,\sigma}$, we will
denote by $\nabla H$ the gradient of $H$ with respect to the scalar
product \eqref{sca.1}, namely
\begin{eqnarray}
\label{e.3}
\nabla H:=\left(\frac{\partial H}{\partial I},\frac{\partial
  H}{\partial \phi} , \frac{\partial H}{\partial p}, \frac{\partial
  H}{\partial q} \right)\ ,
\end{eqnarray}
and by $X_H:=\Jc\nabla H$ the Hamiltonian vector field.

In $\Ph_{s,\sigma}$ consider $r$ Hamiltonian functions $H_\mu^{(l)}$,
$l=1,...,r$ of the form
\begin{equation}
\label{e.4}
H^{(l)}_\mu=I_l+\sum_{j\geq r+1}\Omega_j^{(l)}\frac{p_j^2+q_j^2}{2}+\mu
F^{(l)}_\mu(z)
\end{equation}
When $\mu=0$ the submanifold
\begin{equation}
\label{e.4.2a}
N:=\U\times \T^r\times \left\{0\right\}\times\left\{0\right\}
\end{equation}
is foliated in tori which are invariant under the flow of each
of the Hamiltonians. We are interested in the
persistence of such tori as invariant tori of each of
the Hamiltonians $H^{(l)}_\mu$, when $\mu$ is different from zero.

\begin{remark}
\label{r.sys} Consider $r$ commuting Hamiltonians
$K^{(1)},...,K^{(r)}$, and assume that there exists a torus which is
invariant under the flow of each of the fields
$X_{K^{(l)}}$. Assume also that the fields $X_{K^{(l)}}$ are
linearly independent on such a torus. Then, working as in the proof
of Arnold Liouville theorem, one can prove that, on the torus, the
flow of each of the systems is quasiperiodic.
\emph{Moreover, the torus is invariant under the Hamiltonian flow of
any Hamiltonian which is a linear combination of the $K^{(l)}\null
's$}.
\end{remark}

\begin{remark}
\label{r.red}
In the above situation, assume also that the invariant torus is
elliptic for each flow. If the phase space is finite dimensional, then
it can be shown \cite{Kuk92} that there always exist coordinates in
which the Hamiltonians have the form
\begin{equation}
\label{e.4.a}
K^{(l)}=\sum_{j=1}^{r}\omega_j^{(l)}I_j+\sum_{j\geq
  r+1}\omega_j^{(l)}\frac{p_j^2+q_j^2}{2} +{\rm h.o.t.}\ ,
\end{equation}
where h.o.t. are terms which are of higher order in $I,p,q$ and depend
on $\phi$ also. {Then,
one can construct linear combination of the $K^{(l)}$'s having the form
\eqref{e.4}. So, in this case,} it is equivalent to consider the
Hamiltonians \eqref{e.4.a} {or} the Hamiltonians \eqref{e.4}.  In
\cite{K2} Kuksin studied in detail the case of PDEs, showing that
coordinates in which \eqref{e.4.a} holds exist also in quite general
PDE cases. For this reason here we will directly assume the form
\eqref{e.4}.
\end{remark}

In what follows we will often omit the index $\mu$ of the Hamiltonians
and of the functions $F^{(l)}_\mu$ provided this does not create
ambiguities.

We assume that

{\it
\noindent{(A.1)}
For any $\mu$ small enough one has
\begin{equation}
\label{e.5}
\left\{H^{(l_1)}_\mu;H^{(l_2)}_\mu\right\}=0\ ,\quad \forall l_1,l_2=1,...,r\ .
\end{equation}
}

{\it \noindent(A.2) There exist $d\geq0$ and $s_*\in\R$ such that
$\nabla F^{(l)}_\mu\in C^{\infty}(\V_{s,\sigma},\Ph_{s+d,\sigma})$ for
all $s>s_*$ and for some fixed $\sigma\geq 0$; here
$\V_{s,\sigma}\subset \Ph_{s,\sigma}$ is an open neighborhood of $N$.
}

\begin{remark}
\label{r.smoo}
In (A.2) we are assuming that the nonlinear part of the vector field
is smoothing. Precisely we assume that the nonlinearity allows to gain
$d$ derivatives. This is typical when the system comes from a second
order in time equation.
\end{remark}

We come now to the nonresonance assumption that generalizes to the
case of PDEs Nekhoroshev's nonresonance assumption.

For $\gamma>0$, $\tau\in\R$, and $\bn:=(n_1,...,n_r)\in\Z^r$ consider
the set $\Nc(\gamma,\tau,\bn)$ of the
$(\epsilon_1,...,\epsilon_r)\in\R^r$ such that
\begin{equation}
\label{e.15}
\left|k\sum_{l=1}^{r}n_l-\sum_{l=1}^{r}(n_l+\epsilon_l)
\Omega_j^{(l)}\right| \geq \frac{\gamma}{j^{\tau}}\ ,\quad \forall
k\in\Z\ ,\quad j\geq r+1\ .
\end{equation}

{\it \noindent (A.3) There exist $\bn\in\Z^r$, $\tau\leq d$ and
  $\gamma>0$ such that, for any open set $\Oc\subset\R^r$, whose
  closure contains zero, one has that $\Nc(\gamma,\tau,\bn)\cap \Oc$
  has zero as an accumulation point.  }

\begin{remark}
\label{r.4}
Assumption (A.3) means that, for some $\bn$ the set of the frequencies
which are sufficiently nonresonant accumulates at the origin from any
direction. Moreover we need $\tau\leq d$, with $d$ the amount of
smoothing of the nonlinear part of the vector field.
\end{remark}

Finally we need a non degeneracy condition. In order to state it we
consider the average of the nonlinearity under a suitable periodic
flow. Let $\bn\in\Z^r$ be an integer vector. Define
\begin{equation}
\label{e.8}
F_{\bn}:=\sum_{l=1}^{r}n_lF_\mu^{(l)}\big|_{\mu=0}\ ,
\end{equation}
and
\begin{equation}
\label{e.9}
\langle F_{\bn}\rangle(a,\psi):= \frac1{2\pi}\int_0^{2\pi}F_{\bn}
(a,\bn t+\psi,0,0)dt\ .
\end{equation}

First of all we have the following lemma which will play an important
role in the proof of our main result.

\begin{lemma}
\label{dec}
For any $\bn\in\Z^r$ the function $\langle F_{\bn}\rangle(a,\psi)$ is
independent of $\psi$.
\end{lemma}
\proof The proof is based on a representation formula for the
functions $F_0^{(l)}$ that we are now going to derive.
First remark that assumption (A.1), written explicitly, takes the
form
\begin{eqnarray}
\label{dec.1}
0=\mu\left[ \frac{\partial F^{(l_1)}_\mu}{\partial
\phi_{l_2}}-\frac{\partial F^{(l_2)}_\mu}{\partial \phi_{l_1}}\right]
+\mu\left\{ H^{(l_1)}_{T};F^{(l_2)}_\mu\right\}+\mu\left\{
F^{(l_1)}_\mu;H^{(l_2)}_{T}\right\}
\\
\nonumber
 +\mu^2 \left\{ F_\mu^{(l_1)}; F_\mu
^{(l_2)}\right\}
\end{eqnarray}
where we denoted
$$
H^{(l)}_T:=\sum_{j\geq r+1}\Omega_j^{(l)}\frac{p_j^2+q_j^2}{2}\ .
$$ Restrict the equation \eqref{dec.1} to the manifold $N$.  The terms
containing $H^{(l)}_T$ vanish since they are at least linear in $p$
and $q$. Thus on this submanifold eq. \eqref{dec.1} takes the form
\begin{equation}
\label{dec.2}
\frac{\partial F^{(l_1)}_\mu}{\partial
\phi_{l_2}}-\frac{\partial F^{(l_2)}_\mu}{\partial \phi_{l_1}}=-\mu
\left\{ F_\mu^{(l_1)}; F_\mu^{(l_2)}\right\} \ .
\end{equation}
Evaluating at $\mu=0$ one gets
\begin{equation}
\label{dec.3}
\frac{\partial F^{(l_1)}_0}{\partial
\phi_{l_2}}-\frac{\partial F^{(l_2)}_0}{\partial \phi_{l_1}}=0\ .
\end{equation}
Using the fact that the generators of the first cohomological class of
the torus are $d\phi_l$, $l=1,...,r$, one has that there exists a
function $V$ on $\T^r$ parametrically depending on $I$, and a set of
functions $c_l(I)$, independent of $\phi$, such that
\begin{equation}
\label{dec.4}
F^{(l)}_0(I,\phi)=\frac{\partial V}{\partial \phi_l}+c_l(I)
\end{equation}
Inserting such a representation formula in the definition of $\langle
F_{\bn}\rangle$ one gets that such a function is the sum of two terms, one
of which is independent of $\phi$. The other term contains $V$ and is
proportional to
\begin{equation}
\label{dec.7}
\sum_l\int_0^{2\pi}n_l\frac{\partial V}{\partial \phi_l}(\bn t+\psi,
I)dt =\int_0^{2\pi} \frac{d}{dt}V(\bn t+\psi,I)dt=V(\bn 2\pi+\psi)-V(\psi)=0
\end{equation}
which implies that $\langle
F_{\bn}\rangle(I,\psi) =\sum_ln_lc_l(I)$.\qed

For this reason we will simply write $\langle
F_{\bn}\rangle(a,\psi)\equiv \langle F_{\bn}\rangle (a)$.

{\it \noindent (A.4) There exists $\bn\in \Z^r$ such that (A.3) is
  fulfilled and
\begin{equation}
\label{e.10}
det\left(\frac{\partial^2\langle F_{\bn}\rangle }{\partial
  a^2}(a)\right)\not=0 \ ,\quad \forall a\in\U
\end{equation}
and $\frac{\partial\langle F_{\bn}\rangle }{\partial a}:\U\to\R^r$ is
  a 1-1 map.  }

Fix a set $\U'$ whose closure is contained in $\U$.  Define
\begin{equation}
\label{e.34}
\A(\bn):=\bigcup_{|\mu|<1 } \bigcup_{a\in\U'} \left\{- \mu
  \frac{\partial \langle
F_{\bn}
\rangle }{\partial a}(a)\right\}\subset \R^r
\end{equation}
and assume that this is an open set. Since it has zero as an
accumulation point, it has a nonempty intersection with
$\Nc:=\Nc(\gamma,\tau,\bn)$.  For any $\mu$ consider the set $\A(\bn
)\cap\Nc\cap B_\mu$. Moreover, for any $\epsilon\in \A(\bn
)\cap\Nc\cap B_\mu$ define $a_{*}(\epsilon/\mu)$ as the solution of
\begin{equation}
\label{e.34.b}
\frac{\epsilon}{\mu}=- \frac{\partial \langle F_{\bn} \rangle
  }{\partial a}(a_*)
\end{equation}
We use it to define a reference torus
\begin{equation}
\label{ref.t}
\T_{\epsilon,0}:= \bigcup_{\psi\in\T^r} \left\{(a_*,\psi,0,0)  \right\}\ .
\end{equation}
 Then we have the following theorem.

\begin{theorem}
\label{main} Fix $\sigma$, then for any $s> s_*$ there exist two
constants $\mu_*$ and $C$ such that the following holds.
 For any $|\mu|<\mu_*$ any $\epsilon\in \A\cap\Nc\cap
B_\mu$, there exists a smooth torus $\T_{\epsilon,\mu}\subset
\Ph_{s,\sigma}$ which is invariant under the flow of the Hamiltonian
vector fields of each of the Hamiltonians
$H^{(l)}_\mu$. Moreover one has
\begin{equation}
\label{m.2}
d(\T_{\epsilon,\mu},\T_{\epsilon,0})< C\mu
\end{equation}
where $d(.,.)$ denotes the Hausdorff distance in $\Ph_{s,\sigma}$.
\end{theorem}

\begin{remark}
\label{arnold} As in the first step of the proof of Arnold
Liouville's theorem one has that the flow of each of
the Hamiltonians on the torus $\T_{\epsilon,\mu}$ is the Kronecker
linear flow.
\end{remark}

\begin{remark}
\label{tilde}
From the proof it will also turn out that the flow of
\begin{equation}
\label{ti.1}
\tilde H_\mu:=\sum_l(n_l+\epsilon_l)H_\mu^{(l)}
\end{equation}
on $\T_{\epsilon,\mu}$ is periodic of period $2\pi$ and has trajectories
homotopic to the curve $\phi(t)=\bn t$.
\end{remark}

In all the applications we know one is interested in the dynamics of
one system while the other Hamiltonians are just the generators of the
linear symmetries of the system. In such a case we also compute the
frequencies of the dynamics on the torus.

Precisely assume that there exist $r-1$ independent linear
combinations $K^{(2)},...,K^{(r)}$ of the functions $H^{(l)}_\mu$ wich
are of the form
\begin{equation}
\label{lin}
K^{(l)}=\sum_{j=1}^{r}\omega_j^{(l)}I_j+\sum_{j\geq
  r+1}\omega_j^{(l)}\frac{p_j^2+q_j^2}{2}\ ,
\end{equation}
(no nonlinear part) and have the property that the flow of each
of these Hamiltonians is periodic with period $2\pi$.
Since these functions have no nonlinear part, the validity of
assumption (A.4) implies that $\sum_{l}n_lH^{(l)}_\mu$ and the
functions $K^{(l)}$ are independent. It follows that a system of $r$
independent Hamiltonians in involution is given by $\tilde H_\mu,
K^{(2)},...,K^{(l)}$. It is thus immediate to obtain the following
corollary.

\begin{corollary}
\label{c.main}
Given an arbitrary vector $\alpha=(\alpha_1,...,\alpha_r)\in\R^r$
consider the Hamiltonian $H^{\alpha}_\mu:=\alpha_1 \tilde
H_\mu+\sum_{l=2}^{r}\alpha_l K^{(l)}$. The flow of $H^{\alpha}_\mu$ on
$\T_{\epsilon,\mu}$ is a Kronecker flow with frequency vector
$\alpha$.
\end{corollary}

\section{Proofs}\label{proof}

\begin{remark}
\label{r.13} Fix an integer vector $\bn$ {such that} assumptions
(A.3), (A.4) hold. Make the change of variables
$$
\hat I_1:=\sum_{l}n_lI_l\ ,\quad \hat I_j=I_j\ ,\quad j=2,...,r\ ,
$$ and complete it to a canonical transformation.  Define
\begin{eqnarray}
\label{e.14}
\hat
H^{(1)}&:=&\sum_{l=1}^{r}n_lH^{(l)}
\\
=
\sum_{l=1}^{r} n_l
I_l&+& \sum_{j\geq r+1}\sum_{l=1}^{r}
n_l\Omega^{(l)}_j\frac{p_j^2+q_j^2}{2}+\mu F_n \\ &=&\hat I_1+
\sum_{j\geq r+1} \hat \Omega^{(1)}_j \frac{p_j^2+q_j^2}{2}+\mu \hat
F_1 \ .
\end{eqnarray}
Thus the system of Hamiltonians $\hat H^{(1)}$, $H^{(2)}$,...,
$H^{(r)}$ fulfills the assumptions (A.3) and (A.4)  with $\bn=\be$.
\end{remark}

For this reason the proof will be {carried out} in the case where
(A.3) and (A.4) are fulfilled with $\bn=\be$.  {\it Moreover, from now on
we fix the values of $\gamma$ and $\tau$, and we will simply write
$\Nc$ instead of $\Nc(\gamma,\tau,\be)$}.

In the following we will denote by $\Phi_{\mu,l}^{t_l}$ the time $t_l$
flow of the Hamiltonian vector field $X_{H ^{(l)}_\mu}$. Moreover,
following the standard notation we will use the notation
\begin{equation}
\label{e.15.1}
\Phi^{(t_1,...,t_r)}_\mu:=\Phi_{\mu,1}^{t_1}\circ...\circ
\Phi_{\mu,r}^{t_r}\ .
\end{equation}
$\Phi^t$ is an action of $\R^n$ on $\Ph_{s,\sigma}$.

The first (and main) step of the proof consists in looking for
$r$-tuples $(\epsilon_1,...,\epsilon_r)$ and for initial data $z_0$
such that the function
\begin{equation}
\label{e.16}
z_\mu(t):=\Phi_\mu^{(1+\epsilon_1)t,\epsilon_2t,...,\epsilon_rt}(z_0)
\end{equation}
is periodic of period $2\pi$. Moreover, we will assume that
\eqref{e.16} is contained in a neighborhood of order $\mu$ of a
suitable reference periodic function.

To come to a precise statement consider again $a_*(\epsilon/\mu)$,
namely the solution of \eqref{e.34.b} with $\bn=\be$. Correspondingly
we define the reference periodic function
\begin{equation}
\label{e.34.d}
z_{0,\psi}(t)=(a_*,\be t+\psi,0,0)\subset N\ .
\end{equation}
Then we need to measure the distance between periodic functions.
Thus consider the space
$\Hc_{s,\sigma}:=H^{1}(\T;\R^r\times\R^r\times\ell^2_{s,\sigma}\times
\ell^2_{s,\sigma})$ of the $H^1$ periodic functions of period $2\pi$
{taking} values in the covering space of the
phase space. In such a space we will use the norm
\begin{equation}
\label{e.19.1}
\norma{\zeta}^2_{T,s,\sigma}:=\int_{-\pi}^\pi\norma{\zeta(t)}_{s,\sigma}^2dt +
\int_{-\pi}^\pi\norma{\dot \zeta(t)}_{s,\sigma}^2dt\ .
\end{equation}
Then we consider the space of the functions
of the form
\begin{equation}
\label{e.34.g}
 z(t)=(I(t),\be t+\phi(t),p(t),q(t))=(0,\be t,0,0)+\zeta(t) \ ,
\end{equation}
with $\zeta(t)=(I(t),\phi(t),p(t),q(t))\in \Hc_{s,\sigma}$.
The space of the functions of the form \eqref{e.34.g} will be
symbolically denoted by $\be t+\Hc_{s,\sigma}$. Moreover in $\be
t+\Hc_{s,\sigma}$ we will use the topology induced by the distance
introduced by the norm of $\Hc_{s,\sigma}$.

\begin{lemma}
\label{l.main} Fix $\sigma$, then, for any $s> s_*$ there exist two
constants $\mu_*$ and $C$ such that the following holds.
 For any $|\mu|<\mu_*$, any $\epsilon\in
\A(\be)\cap\Nc\cap B_\mu$, and any $\psi\in\T^r$, there exists a
\emph{unique} $2\pi$ periodic function $z_{\mu,\psi}^\epsilon(t)\in
\be t+\Hc_{s,\sigma}$ of the form \eqref{e.16} fulfilling
\begin{equation}
\label{e.34g}
\norma{z_{0,\psi}-z_{\mu,\psi}^\epsilon}_{T,s,\sigma }\leq C\mu
\end{equation}
and
$$
z_{\mu,\psi}^\epsilon(0)=(I(0),\psi,p(0),q(0))
$$ i.e. the initial datum for the component on $\T^r$ is exactly
$\psi$.  Moreover the map $(\mu,\psi)\mapsto z_{\mu,\psi}^\epsilon\in
\be t+\Hc_{s,\sigma}$ is $C^\infty$.
\end{lemma}

The proof of this lemma consists of several steps and is contained in
the next subsection

\subsection{Proof of lemma \ref{l.main}}\label{ss.1}

A function of the form \eqref{e.16} is a solution of the
Hamiltonian system
\begin{eqnarray}
\label{e.17}
\tilde H:=(1+\epsilon_1)I_1+\epsilon_2I_2+...+\epsilon_rI_r+
\sum_{j\geq r+1}\tilde \Omega_j(\epsilon)\frac{p_j^2+q_j^2}{2}+\mu
\tilde F_{\mu}(z)
\end{eqnarray}
where
\begin{eqnarray*}
\tilde \Omega_j(\epsilon):=(1+\epsilon_1)\Omega_j^{(1)}+
\sum_{l=2}^{r}\epsilon_l\Omega^{(l)}_j \ ,\quad \tilde F_{\mu}
:=(1+\epsilon_1)F^{(1)}_\mu+\sum_{l=2}^{r}\epsilon_lF^{(l)}_\mu\ ;
\end{eqnarray*}
so we will actually look for a periodic solution of such a
Hamiltonian system.

As usual when dealing with systems having some integrals of motion the
situation is quite delicate. Thus, instead of looking for periodic
solutions of $\tilde H$ we look for constants $\beta_l$ and for
periodic solutions of the system
\begin{equation}
\label{e.18}
\dot z=X_{\tilde H}(z)+\sum_{l=1}^{r}\beta_l\nabla\Ic^{(l)}(z)\ ,
\end{equation}
where $\Ic^{(l)}(z):=I_l$.  This is possible in view of the following
lemma.

\begin{lemma}\label{3.2}
Consider a function $z \in \be t+\Hc_{s,\sigma}$, $s>s_*$. Then there
exists a $\mu_*>0$ such that, if $z(t)\in\V_{s,\sigma}$ (cf. (A.2)),
$\forall t$, and if $|\mu|<\mu_*$, then $(z(t),\beta)$ is a solution
of the system \eqref{e.18} if and only if $\beta=0$ and $z(t)$ is a
periodic solution of the system
$$
\dot z=X_{\tilde H}(z)\ .
$$
\end{lemma}
\proof Assume that $z(t)$ is a periodic solution
of \eqref{e.18}, and apply $d H^{(l_1)}(z(t))$ to \eqref{e.18}. One
has
\begin{eqnarray*}
 d H^{(l_1)}(z(t))\dot z=d H^{(l_1)}X_{\tilde H}+\sum_{l=1}^{r}\beta_ld
H^{(l_1)} \nabla \Ic^{(l)}(z) \\ = \left\{H^{(l_1)} ; \tilde H \right\}
+\sum_{l=1}^{r}\beta_l \langle \nabla H^{(l_1)};\nabla \Ic^{(l)}\rangle
\\
=\sum_{l=1}^{r}\beta_l \langle \nabla H^{(l_1)};\nabla
\Ic^{(l)}\rangle
=\sum_{l=1}^{r}\beta_l\frac{\partial H^{(l_1)}}{\partial I_l}
\\
=
\sum_{l=1}^{r} \beta_l\left(\delta_{l,l_1}+\mu\frac{\partial
  F^{(l_1)}}{\partial I_l}\right)\ .
\end{eqnarray*}Take the integral over $2\pi$ of such an equation, and
{notice} that the left hand side is the time
derivative of $H^{(l_1)}(z(t))$, which by assumption is a periodic
function, thus we get
\begin{eqnarray}
\label{e.181}
0=\sum_{l=1}^{r}\beta_l\left[ \int_{0}^{2\pi}
  \left(\delta_{l,l_1}+\mu\frac{\partial F^{(l_1)}}{\partial
  I_l}(z(t))\right) dt \right]\ .
\end{eqnarray}
Then, for $\mu$ small enough the square bracket is an invertible
matrix, and therefore \eqref{e.181} implies $\beta=0$. \qed

We look now for $\beta$, $\epsilon$ and a periodic solution
$z(t)\in\be t+\Hc_{s,\sigma}$ of \eqref{e.18}. Write $z(t)=(0,\be
t,0,0)+\zeta(t)$ with $\zeta\equiv
(I(t),\phi(t),p(t),q(t))\in\Hc_{s,\sigma}$; substituting in
\eqref{e.18} one gets the system:
\begin{equation}
\label{e.304} L_{\epsilon}\zeta =\Theta (\zeta,\epsilon,\beta,\mu),
\end{equation}
where the operator $L_\epsilon$ is defined by
\begin{equation}
\label{e.25a}
L_{\epsilon}\left(
\begin{matrix}
I_j \\ \phi_j \\ p_j \\ q_j
\end{matrix}
\right)
=
\left[
\begin{matrix}
\frac{d}{dt} &0&0&0
\\
0& \frac{d}{dt} &0&0
\\
0&0&\frac{d}{dt} & \tilde
\Omega_j(\epsilon)
\\
0&0&- \tilde
\Omega_j(\epsilon) &\frac{d}{dt}
\end{matrix}
\right]
\left(
\begin{matrix}
I_j \\ \phi_j \\ p_j \\ q_j
\end{matrix}
\right)
\end{equation}
and $\Theta$ is a nonlinear operator
\begin{eqnarray}
\label{e.191}
\Theta:\Hc_{s,\sigma}\times \R^r\times\R^r\times \R &\to&
\Hc_{s+d,\sigma} \\ (\zeta,\epsilon,\beta,\mu) &\mapsto& \zeta'\equiv
\Theta(\zeta,\epsilon,\beta,\mu )
\end{eqnarray}
with $\zeta'(t)=( I'(t),\phi'(t),p'(t), q'(t))$ defined by
\begin{eqnarray}
\label{e.19.b}
I'_j(t)=-\mu\frac{\partial \tilde
  F}{\partial\phi_j}(\zeta(t))+\beta_j\ ,
\qquad\qquad
\phi'_j(t)=\epsilon_j+ \mu \frac{\partial \tilde
  F}{\partial I_j} (\zeta(t))
\\
\label{e.19.c}
p_j'(t)
= -\mu \frac{\partial \tilde
  F}{\partial q_j}(\zeta(t))
\ ,\qquad
q'_j = \mu \frac{\partial \tilde F}{\partial p_j}(\zeta(t))\ .
\end{eqnarray}
\begin{remark}
\label{r.30}
Since $\Hc_{s,\sigma}$ is an algebra, by assumption
(A.2) the operator $\Theta$ is a smooth map from $\Hc_{s,\sigma}\to
\Hc_{s+d,\sigma}$ for all $s>s_*$.
\end{remark}
We solve the system \eqref{e.304} by using Lyapunof-Schmidt
decomposition. It is easy to see that the kernel of
$L_0:=L_\epsilon\big|_{\epsilon=0}$ is given by the space of the
constant functions of the form $(a,\psi,0,0)$, while its range is the
space $R_{s,\sigma}$ of the functions $w (t)\equiv
(J(t),\theta(t),p(t),q(t))\in \Hc_{s,\sigma} $ with $J(t),\theta(t)$
having zero average.

We project our system of equations on the two subspaces $R_{s,\sigma}$
and $K$. Correspondingly we also decompose the unknown as
$\zeta(t)=(a,\psi,0,0)+ w (t)$, $w\in R_{s,\sigma}$.

Denote by $P$ the projector on $R_{s,\sigma}$, and by $Q$ the
projector on $K$. Explicitly $Q$ acts by putting equal to zero the
$p$ and $q$ components of $\zeta$ and by taking the average of the
first two components.  Thus the system \eqref{e.304} decomposes into
the range equation
\begin{equation}
\label{e.305}
L_{\epsilon}w  =P\Theta (w+(a,\psi,0,0),\epsilon,\beta,\mu)\ ,
\end{equation}
and the kernel system
\begin{eqnarray}
\label{e.22}
0&=&-\mu\left\langle\frac{\partial \tilde
  F}{\partial\phi_j}\big(w(t)+(a,\be t+\psi,0,0)\big) \right\rangle+\beta_j
\ ,
\\
\label{e.22.1}
0&=&\epsilon_j+\mu \left\langle\frac{\partial \tilde F}{\partial I_j}
  \big(w(t)+(a,\be t+\psi,0,0)\big)\right\rangle
\end{eqnarray}
Where we denoted by $\langle . \rangle$ the operation of taking the
average, for example
\begin{equation}
\label{ave}
\langle J_j\rangle=\frac{1}{2\pi}\int_0^{2\pi}J_j(t)dt\ .
\end{equation}
First fix $\epsilon\in \Nc$, an open set $\U''\supset \U'$, with
closure contained in $\U$. We fix $a\in \U''$, $\beta$ in a
neighborhood of zero and $\psi\in\T^r$ and we solve the system
\eqref{e.305} getting a solution $w_\epsilon(\mu,a,\psi)$ which turns
out to be independent of $\beta$. Then we substitute in the system
\eqref{e.22}-\eqref{e.22.1} and we solve it.

\vskip20pt Consider \eqref{e.305}; we have that under the assumption
(A.3) the operator $L_{\epsilon}$ admits an inverse which is bounded
as an operator from $ R_{s+\tau,\sigma}$ to $ R_{s,\sigma}$. Precisely
the following lemma holds.

\begin{lemma}
\label{l.3}
Let $w \in R_{s+\tau,\sigma}$, fix $\epsilon\in \Nc$
then there exists a unique $w '\in R_{s,\sigma}$ such that
\begin{equation}
\label{e.25}
L_{\epsilon}w '=w \ ,
\end{equation}
moreover one has
\begin{equation}
\label{e.26}
\norma{w '}_{{T,s,\sigma}}\leq \frac{1}{\gamma}
\norma{w }_{{T,s+\tau,\sigma}}
\end{equation}
\end{lemma}
\proof {Notice} that the operator $L_\epsilon$ is block diagonal.
Denote $w =(J_j,\theta_j,p_j,q_j)$ and similarly for $w '$, expand
each of the components of $w $ in time Fourier series, i.e. write
\begin{equation}
\label{e.27}
J_j(t)=\frac{1}{\sqrt{2\pi}}\sum_{k}J_{jk}e^{ikt}
\end{equation}
and similarly for all the other components. One has that {the} system \eqref{e.25} is equivalent to the system
\begin{eqnarray}
\label{e.28}
ik J'_{jk}:={J_{jk}}\ ,\quad ik \theta'_{jk}:={\theta_{jk}}\ ,\qquad
j=1,...,r\ ,\ k\not=0 \\
\label{e.29}
\left[\begin{matrix}
ik &\tilde \Omega_j
\\
-\tilde \Omega_j & ik
\end{matrix}\right]\left(
\begin{matrix}
p'_{jk}\\q'_{jk}
\end{matrix}\right)  =\left(
\begin{matrix}
p_{jk}\\q_{jk}
\end{matrix}\right)\ ,\qquad j\geq r+1\ ,\ k\in\Z \ .
\end{eqnarray}
{Notice} also that in terms of the Fourier
variables the norm takes the form
\begin{eqnarray}
\nonumber
\norma{ w}_{T,s,\sigma}^2=\sum_{k}(1+k^2)\left[\sum_{j=1}^{r} \left(
  \left|J_{jk}\right|^2 + \left|\theta_{jk}\right|^2 \right)
\right.
\\
\label{no.e}
\left.+\sum_{j\geq r+1} \left( \left|p_{jk}\right|^2 +
  \left|q_{jk}\right|^2 \right)[j]^{2s}e^{2\sigma|k|} \right]
\end{eqnarray}
It is immediate to solve the two equations \eqref{e.28}. To solve
\eqref{e.29} one can simply observe that the matrix defining such a
system has eigenvalues and eigenvectors given by
\begin{equation}
\label{e.30}
\lambda_{kj}=i\left( k\pm\tilde \Omega_j(\epsilon)\right)\ ,\quad
\frac{1}{\sqrt2}\left(
\begin{matrix}
1 \\ \pm i
\end{matrix}
\right)\ .
\end{equation}
In particular the two eigenvectors are orthogonal {to} each
other and independent of $k$ and $j$. With this information it is
easy to construct and estimate the solution of \eqref{e.29} getting
\begin{equation}
\label{e.31}
\norma{(p',q')}_{T,s,\sigma}\leq \norma{(p,q)}_{T,s+\tau,\sigma}
\sup_{jk}\left[\frac{1}{j^\tau\left|i(k\pm\tilde \Omega_j(\epsilon))
\right|} \right]
\end{equation}
If $\epsilon\in\Nc$ then the argument of the supremum is
estimated by
$$ \frac{1}{j^\tau\left|k\pm\tilde \Omega_j(\epsilon) \right|}\leq
\frac{j^{\tau}}{j^\tau \gamma}=\frac{1}{\gamma}\ .
$$
Adding the trivial estimate of $J'$ and of $\theta'$ one gets the
thesis.\qed

According to the previous lemma one can define
$L_{\epsilon}^{-1}$ which is bounded as an operator from
$R_{s+\tau,\sigma}$ to $R_{s,\sigma}$. It follows that one can rewrite
the system \eqref{e.25} as a fixed point problem, namely
\begin{equation}
\label{e.32}
w =
L_{\epsilon}^{-1}P\Theta(w+(a,\psi,0,0),\epsilon,\beta,\mu)\equiv
\mu\F_{\epsilon} (w ,\mu,a,\psi)
\end{equation}
where $\F$, is a map which, for fixed $\epsilon\in\Nc$, is a smooth
map from $R_{s,\sigma}\times (-\mu_*,\mu_*)\times\U''\times \T^r$ to
$R_{s,\sigma}$. We factorized a $\mu$ in front of $\F$, since it
explicitly appears in the form of $P\Theta$ and eliminated $\beta$
since $P\Theta$ is indepndent of it.

So one can apply the implicit function {theorem} in order to get
the following corollary.

\begin{corollary}
\label{c.1}
Fix $s>s_*$ and $\sigma$, then, for any $\epsilon\in \Nc$ there exists
a function $w (\mu,a,\psi )$ which solves the system
\eqref{e.32}. Moreover it is $C^{\infty}$ as a map from
$(-\mu_*,\mu_*)\times\U''\times \T^r $ to $R_{s,\sigma}$, and fulfills
the inequality
\begin{equation}
\label{e.33}
\norma{w (\mu,a,\psi)}_{{T,s,\sigma}}\leq \mu C
\end{equation}
for all $(\mu,a,\psi)$ in the considered domain.
\end{corollary}

Then we substitute such a solution in \eqref{e.22}-\eqref{e.22.1} and
we solve such a system. Since $w $ does not depend on $\beta$, the
solution of \eqref{e.22} is immediate. The solution of \eqref{e.22.1}
is slightly more difficult since $w $ does not depend smoothly on
$\epsilon$ and moreover it is defined only for $\epsilon$ in the set
$\Nc$. For this reason we proceed as in \cite{bam00}, i.e. we fix
$\epsilon\in\Nc$ and we look for a vector $a=a_\epsilon(\mu,\psi)$
which fulfills such an equation.

We
have the following lemma.

\begin{lemma}
\label{sol.Q} There exists a positive $\mu_*$ such that for any
$|\mu|<\mu_*$ the following holds.  For any
$\epsilon\in\A\cap\Nc\cap B_{\mu}$ there exists a unique smooth
function $\alpha_{\frac{\epsilon}{\mu}}(\mu,\psi)$ such that
\begin{equation}
\label{e.35}
a=a_*+\mu\alpha_{\frac{\epsilon}{\mu}}(\mu,\psi)
\end{equation}
solves the equation \eqref{e.22.1}.
\end{lemma}
\begin{remark}
\label{r.er}
The function $\alpha$ is not smooth in $\epsilon/\mu$.
\end{remark}

\proof First {notice} that from \eqref{e.33}
one has
\begin{equation*}
\begin{array}{c}
\nonumber \left\langle \!\frac{\partial\tilde F}{\partial
I_j}\left((a,\be t+\psi,0,0)+w(\mu,a,\psi,t)\right)\! \right\rangle
\!=\!\left\langle\! \frac{\partial\tilde F}{\partial I_j}(a,\be
t+\psi,0,0)\! \right\rangle\!+\!\mu \G(a,\psi,\mu)
\end{array}
\end{equation*}
with  $\G(a,\psi,\mu) $ a suitable smooth function. Thus the above
expression is equal to
\begin{eqnarray}
\frac{1}{2\pi}\!\int_0^{2\pi}\!\!\frac{\partial\tilde F}{\partial
  I_j}(a,\be t+\psi,0,0)dt+\mu \G(a,\psi,\mu)
\!=\!\frac{\partial\langle \tilde
  F\rangle }{\partial
  a_j}(a)+\mu \G(a,\psi,\mu)\,\,\,\,\,\,
\label{e.322}
\end{eqnarray}
where of course
$$
\langle\tilde
  F\rangle (a)=\frac{1}{2\pi} \int_0^{2\pi}\tilde F(a,\be
    t+\psi,0,0)dt\ .
$$ 
Write $a=a_*+\mu\alpha$, thus we have to solve
$$ 0=\frac{\epsilon}{\mu}+\frac{\partial\langle \tilde F\rangle
}{\partial a_j}(a_*+\mu\alpha)+\mu
\G(a_*+\mu\alpha,\psi,\mu)=:\Lc(\rho,\mu,\alpha,\psi)
$$
for a fixed value of the parameters $\rho:=\frac{\epsilon}{\mu}$.
Due to the definition of $a_*$ we have
$$
\norma{\Lc(\rho,\mu,\alpha,\psi)}\leq C\mu\ .
$$ Moreover one has that $\frac{\partial \Lc}{\partial \alpha}\simeq
\frac{\partial^2 \langle\tilde F \rangle}{\partial \alpha^2} $ is
invertible with inverse having norm independent of $\mu$. Then the
implicit function theorem ensures the result\footnote{actually one
needs a small variant of the standard implicit function theorem, since
we do not know any solution of the equation $\Lc=0$. However the
standard proof based on the contraction mapping theorem allows one to
prove that the implicit function exists also in this case.}  \qed

\noindent
{\it End of the proof of lemma \ref{l.main}.} Thanks to lemma
\ref{sol.Q} we have constructed the unique periodic function of the
form \eqref{e.16} belonging to $\be t+\Hc_{s,\sigma} $ which is $\mu$
close to $z_{0,\psi}$ and such that {\it the average of the angular
component $\phi$ is equal to $\psi$}. In order to get lemma
\ref{l.main} it is enough to show that there is a 1-1 correspondence
between the average and the initial datum for the angular
component. This is an immediate consequence of the trivial relation
\begin{equation}
\label{e.306}
\phi(0)=\psi+\theta(0)\ ,
\end{equation}
and of the fact that $\theta(0)$ is a smooth function of $\psi$ and
has size
smaller than $\mu$. Thus from the implicit function theorem one can
calculate the average $\psi$ as a unique smooth function of
$\phi(0)$. This concludes the proof. \qed

\subsection{End of the proof of theorem \ref{main}.}

Fix $|\mu|<\mu_*$ and $\epsilon\in\A\cap\Nc\cap B_{\mu}$.  Consider
the torus
\begin{equation}
\label{g.1}
\T_{\epsilon,\mu}:= \bigcup_{\psi\in\T^r}\{ z_{\mu,\psi}^\epsilon(0) \}
\end{equation}
where $z_{\mu,\psi}^\epsilon(t)$ is the periodic solution of $\tilde
H$ constructed in lemma \ref{l.main}. Moreover one has an important
uniqueness property ensuring that there are no other periodic
solutions of $\tilde H$ close to such a torus. We exploit now such a
uniqueness property in order to prove the following lemma which is the
last step for the proof of theorem \ref{main}.

\begin{lemma}
\label{geo} The torus $\T_{\epsilon,\mu}$ is invariant under the
flow of the Hamiltonian vector fields of each of the
functions $H^{(l)}_\mu$.
\end{lemma}
\proof Fix $z\in \T_{\epsilon,\mu}$, and a neighborhood $\U_0\subset
\R^r$ of zero. Consider the set
\begin{equation}
\label{g.2}
\Sc:=\bigcup_{\varphi\in \U_0} \Phi_\mu^{\varphi}(z)\ .
\end{equation}
Since $\Phi^\varphi_\mu$ commutes with the flow of $\tilde H$, all
the points of $\Sc$ give rise to periodic orbits of $\tilde H$.
Moreover $\Sc$ and $\T_{\epsilon,\mu}$ have at least $z$ (and the
periodic orbit it generates) in common. But there are no initial
data outside $\T_{\epsilon,\mu}$ and close to it, which give rise to
periodic solutions of $\tilde H$ (uniqueness property of the torus).
It follows that $\Sc\subset \T_{\epsilon,\mu}$. So this property
implies the thesis. \qed

\section{Applications}\label{Appli}

\subsection{Two dimensional tori in a nonlinear wave equation}
\label{nlw}

Consider the nonlinear wave equation
\begin{equation}
\label{nlw.1}
u_{tt}-u_{xx}+m u=f(u)\ ,\quad x\in\T\equiv\frac{\R}{2\pi\Z}\ .
\end{equation}
where $f$ is a smooth function having a zero of order higher than 1 at
the origin. For simplicity we will consider explicitly only the case
$f(u)=-u^3$ We will construct small amplitude invariant tori of
dimension 2 exploiting the invariance of the equation under
translations.

The equation \eqref{nlw.1} is Hamiltonian with Hamiltonian function
\begin{equation}
\label{nlw.3}
K^{(1)}(U,u)=\int_{-\pi}^\pi\left[\frac{U^2+u_{x}^2+mu^2}{2}
+\frac{u^4}{4}\right] dx
\end{equation}
where $U\equiv \dot u$ is the momentum conjugated to $u$. Translations
in $\T$ are generated by the Hamiltonian
\begin{equation}
\label{nlw.4}
K^{(2)}(U,u)=\int_{-\pi}^\pi U_xudx
\end{equation}
which commutes with $K^{(1)}$. To fit the scheme of the previous
section we have to introduce suitable variables. First introduce the
Fourier basis $\hat e_j(x)$ by
\begin{equation}
\label{nlw.5}
\hat e_j(x)=\left\{
\begin{matrix}
\frac{\cos(jx)}{\sqrt\pi} & {\rm if}\ j>0
\\
\frac{\sin(-jx)}{\sqrt\pi} & {\rm if}\ j<0
\\
\frac{1}{\sqrt{2\pi}} & {\rm if}\ j=0
\end{matrix}
\right.\ ,
\end{equation}
{Notice} that $\partial_x\hat e_j=-j\hat
e_{-j}$ and introduce the variables $u_j$, $U_j$ by
\begin{equation}
\label{nlw.6}
u(x)=\sum_{j\in\Z}\frac{u_j}{\sqrt{\omega_j}}\hat e_j(x)\ ,\quad
U(x)=\sum_{j\in\Z} \sqrt{\omega_j}\ U_j\hat e_j(x)\ ,
\end{equation}
where
\begin{equation}
\label{nlw.7}
\omega_j=\sqrt{j^2+m}\ .
\end{equation}
Then we introduce the variables $p_j,q_j$ by
\begin{eqnarray}
\label{nlw.8}
u_j=\frac{p_j+q_{-j}}{\sqrt 2}\ ,\quad U_j=\frac{p_{-j}-q_{j}}{\sqrt
  2}\ ,\quad j\not=0\ ,
\end{eqnarray}
and $p_0=U_0$, $q_0=u_0$.  It is easy to check that these are
canonical variables and that the Hamiltonians take the form
\begin{eqnarray}
\label{nlw.9}
K^{(1)}&=&\sum_{j\in\Z}\omega_j\frac{p_j^2+q_j^2}{2}+F
\\
K^{(2)}&=&\sum_{j\in\Z}j\frac{p_j^2+q_j^2}{2}
\end{eqnarray}
where
\begin{equation}
\label{nlw.10}
F:=\int_{-\pi}^\pi\frac{u^4}{4} dx
\end{equation}
has to be thought as a function of the variables $p,q$.

Then we rescale the variables introducing the small parameter $\mu$ by
defining $p=\sqrt\mu \tilde p$, $q=\sqrt\mu \tilde q$, we rewrite the
Hamiltonians in terms of such variables, {\it and omit tildes}.

We fix two modes of the linear systems and we continue to the
nonlinear system the family of invariant tori formed by the
superposition of the linear oscillations involving only such two
modes.  Just to be determined (and because the computations turn out
to be simple) we take the two modes with indexes 1 and -1.

Finally we introduce the action angle variables for the two
considered modes by defining the actions
$$ I_{1}=\frac{p_1^2+q_1^2}{2}\ ,\quad
I_{-1}=\frac{p_{-1}^2+q_{-1}^2}{2}\ ,
$$ and the corresponding angles. So the two Hamiltonians $K^{(1)},
K^{(2)}$ turn out to have the form \eqref{e.4.a}. A simple computation
shows that the two functions
\begin{equation}
\label{nlw.12}
H^{(1)}=\frac{K^{(1)}+\omega_1 K^{(2)}}{2\omega_1}\ ,\quad
H^{(2)}=\frac{K^{(1)}-\omega_1 K^{(2)}}{2\omega_1}\ ,
\end{equation}
have the structure \eqref{e.4}, suitable for the application of theorem
\ref{main}. In particular the frequencies are
\begin{equation}
\label{nlw.13}
\Omega^{(1)}_{\pm j}=\frac{\omega_j\pm j\omega_1}{2\omega_1}\ ,\quad
\Omega^{(2)}_{\pm j}= \frac{\omega_j\mp j\omega_1}{2\omega_1}\ ,\quad
j\geq0\ ,\quad j\not=1\ .
\end{equation}
The nonlinearities are $F^{(1)}=F^{(2)}=F/2\omega_1$. We fix the
domain $\U$ to be the open set $(0,+\infty)\times (0,+\infty)$.  We
fix the indexes of the phase space to be some arbitrary $\sigma>0$
and some $s> 1/2$. Then, by the Sobolev embedding theorem the
nonlinearity fulfills assumption (A.2) with $d=1$. We come to the
nonresonance assumption (A.3). It turns out that it is convenient to
choose $n_1=n_2=1$. We have the following
\begin{lemma}
\label{l.nlw.1}
There exists an uncountable dense subset $\Sc$ of $(-1,\frac{4}{3})$,
with zero measure, such that, if $m\in\Sc$ then property (A.3) holds
with $\tau=1$ and a positive $\gamma=\gamma(m)$. Fix $m\in\Sc$, denote
by $\Nc_0$ the intersection of the set $\Nc$ with a given neighborhood
of the origin. Then $\Nc_0$ is uncountable; as
$(\epsilon_1,\epsilon_2)\equiv\epsilon$ varies in $\Nc_0$, the
quantity $\epsilon_1-\epsilon_2$ assumes uncountably many values, but
also rational values.
\end{lemma}
\begin{remark}
\label{rem.nlw.1}
The quantity $\epsilon_1-\epsilon_2$ is important since it will be the
ratio between the two frequencies of motion on the invariant torus.
\end{remark}
\proof One has, for $j\geq 0$,
\begin{eqnarray*}
\tilde \Omega_{\pm j}:=(1+\epsilon_1)\Omega^{(1)}_{\pm
j}+(1+\epsilon_2)\Omega_{\pm j}^{(2)}
=\frac{\omega_j}{\omega_1}\left(1+\frac{\epsilon_1+\epsilon_2}{2}\right)
\pm \frac{\epsilon_1-\epsilon_2}{2}j \\
=\left(\frac{1}{\omega_1}\left(
1+\frac{\epsilon_1+\epsilon_2}{2}\right)\pm
\frac{\epsilon_1-\epsilon_2}{2} \right) j+\frac{m}{\omega_j+j}
\end{eqnarray*}
Consider first the case $\epsilon_1=\epsilon_2=0$. We construct a set
$\Sc_0$ such that, if $m\in\Sc_0$ then $\tilde \Omega_{j}\equiv
\frac{\omega_j}{\omega_1}$ fulfills
\begin{equation}
\label{nlw.16}
\left|k-\tilde\Omega_j\right|\geq\frac{\gamma}{|j|}\ ,\quad \forall
k,j, \quad {\rm with}\ |k|,|j|\geq J_*\ ,\gamma=\frac{1}{6}
\end{equation}
and $J_*$ a large number. We follow closely the construction of
\cite{BP02}, proof of lemma 3.1.  Fix $\bar m\in(-1,4/3)$ and
$\delta>0$; we construct uncountably many values of $m$ which are
$\delta$ close to $\bar m$ and such that \eqref{nlw.16} holds.
Consider $\omega_1(m):=\sqrt{1+m}$, and $\nu(m):=[\omega_1(m)]^{-1}$,
$\bar \nu:=\nu(\bar m)$. {Notice} that there is 1-1 correspondence
between $m$ and $\nu$. Consider the continued fraction expansion of
$\bar \nu$. Fix some large integer $Q$ and consider the set $\Sc_\nu$
of the $\nu$'s obtained by substituting the terms of index larger than
$Q$ in the continued fraction expansion of $\bar \nu$, with an
infinite sequence of $0,1$. If $Q$ is large enough then $\Sc_\nu$ is
contained in an arbitrarily small ball centered at $\bar \nu$. By
standard continued fractions theory one has that if
$\nu\in\Sc_\nu$, then
$$
|k-\nu j|\geq \frac{\gamma'}{j}\quad \forall j>J_*\ ,\quad \forall
k\in\Z \quad \gamma'=\frac{1}{3}\ ,
$$ where $J_*$ depends on $Q$.  Denote by $\Sc_0$ the preimage of
$\Sc_\nu$ under the map $m\mapsto \nu(m)$, and assume that $Q$ is so
large that $\Sc_0$ is contained in a ball of radius $\delta$ centered
at $\bar m$. Thus, taking into account also the quantity
$m/(\omega_j+j)$, we have that, for $m\in\Sc_0$,
$$ |k-\tilde \Omega_j|\geq \frac{\gamma}{j}\ ,\quad \forall j>J_*\
,\quad \forall k\ .
$$ Finally we take out of $\Sc_0$ the set of the $m$'s
such that $k-\tilde \Omega_j=0$ for some $j\leq J_*$, $k\in\Z$. This is a
finite set since $\nu$ is an analytic function of $m$. Thus the
remaining set is contained in a ball of radius $\delta$ centered at
$\bar m$, and it is uncountable. The union of such sets over all the
points $\bar m$ is by definition the set $\Sc$.

We discuss now the case $\epsilon\not=0$. Fix $m\in\Sc$, denote
\begin{eqnarray*}
\epsilon_+:=\frac{\epsilon_1+\epsilon_2}{2}\ ,\quad
\epsilon_-:=\frac{\epsilon_1-\epsilon_2}{2}\ , \\ \nu_+=
\nu(1+\epsilon_+)+\epsilon_-\ ,\quad \nu_-=
\nu(1+\epsilon_+)-\epsilon_-\ ,
\end{eqnarray*}
and {notice} that
$$ \tilde \Omega_{\pm j}=j\nu_{\pm}+\frac{m}{\omega_j+j}\ ,\quad j\geq
0\ .
$$ The maps $(\epsilon_1,\epsilon_2)\mapsto
(\epsilon_+,\epsilon_-)\mapsto (\nu_+,\nu_-)$ are 1-1. Consider the
set $\Nc_1$ of the couples $\epsilon_1,\epsilon_2$ such that $\nu_+$
and $\nu_-$ have a continued fraction expansion $[a_1,a_2,...]$ with
$a_j\leq 1$ for $j$ large enough. Since the couples of $(\nu_+,\nu_-)$
belonging to an arbitrary neighborhood of $(\nu,\nu)$ are uncountable,
the same property holds for $\Nc_1$.  The result on the existence of
rational values of $\epsilon_-=\epsilon_1-\epsilon_2$ is obtained by
fixing $\epsilon_+$ in such a way that $\nu(1+\epsilon_+)$ is of constant
type, and by {noticing} that both $\nu_+$ and $\nu_-$ remain of
constant type when $\epsilon_-$ is a rational number. Thus in $\Nc_1$
one has
$$ |k-\tilde \Omega_j(\epsilon)|\geq \frac{\gamma}{j}\quad \forall
j>J_*\ ,\quad \forall k\ .
$$ To deal with the smaller values of $j$ consider the (finite) set of
functions
$$
j(k-\tilde \Omega_j(\epsilon)) \quad
j\leq J_*\ ,\quad k\leq K_*
$$ they are different from zero at $\epsilon=0$, so there exists a
neighborhood of 0 such that they remain different from zero. The set
$\Nc$ is the intersection of $\Nc_1$ with such a neighborhood of
zero.\qed

Finally we have to prove that the non degeneracy assumption (A.4)
holds. To this end we have to compute the function $F_{\bn}$ and its
average.

First of all we need to compute the flow over which we have to
average. This is easily done by going back through the changes of
variables we did. Thus one gets
\begin{eqnarray*}
u(a,\bn t,0,0)=\sqrt{a_{-1}}\frac{\cos(x-t)}{\sqrt
  \pi}+\sqrt{a_{1}}\frac{\sin (x+t)}{\sqrt \pi}
\end{eqnarray*}
Inserting in $F$, in $F_{\bn}$ and computing $\langle F_\bn\rangle$ one gets
that such a function is
proportional (through a novanishing constant) to
\begin{equation}
\label{nlw.33} a_1^2+{4}a_1a_{-1}+a_{-1}^2
\end{equation}
which is a nondegenerate quadratic form, and thus property (A.4)
holds.

Thus we can apply the general theory and get quasiperiodic solutions
in the nonlinear wave equation. In order to get a slightly more global
description it is convenient to scale back the variables eliminating
the parameter $\mu$. We formulate the result for the nonlinear wave
equation in terms of the original amplitudes $A=\mu a$.

Consider the set
\begin{equation}
\label{nlw.161}
\A_*:=\bigcup_{A\in (0,\infty)\times (0,\infty) }\left\{-\frac{\partial
\langle F_{(1,1)}\rangle}{\partial a}(A) \right\}
\end{equation}
and, for $\epsilon\in \A_*$, denote by $A_*(\epsilon)$ the unique
solution of
$$
\epsilon =-\frac{\partial
\langle F_{(1,1)}\rangle}{\partial a}(A) \ ,
$$
and by $\T_{\epsilon,*}$ the torus $p=q=0, I=A_*, \phi\in\T^2$.
\begin{theorem}
\label{nlw.teo}
There exists a positive $\mu_*$ with the following property: for any
$\epsilon\in \Nc\cap B_{\mu_*}\cap \A_*$ there exists a unique invariant torus
$\T_{\epsilon}$ such that
\begin{itemize}
\item[1.] the flow on $\T_{\epsilon}$ has frequencies
\begin{equation}
\label{fre.nlw}
\left(\frac{\omega_1}{1+\frac{\epsilon_1+\epsilon_2}{2}},\ -
\frac{\omega_1}{1+\frac{\epsilon_1+\epsilon_2}{2}}
\frac{\epsilon_1-\epsilon_2}{2}\right)
\end{equation}
\item[2.]  one has
\begin{equation}
\label{nlw.331}
d(\T_{\epsilon},\T_{\epsilon,*})\leq C \norma{A_*(\epsilon)}^2
\end{equation}
\end{itemize}
\end{theorem}
\proof We apply theorem \ref{main}. To this end we must use the
variables $a$ and the parameter $\mu$. For fixed $\epsilon$ define
$\mu:=\norma{A_*(\epsilon)}$, and fix the set $\U'$ to be the
spherical shell $\norma a=1$. Then using the notation of
sect. \ref{stat} one has $a_*(\epsilon/\mu)=A_*(\epsilon)/\mu$.  Then
remark that in this case one has
$$
K^{(1)}=\frac{\omega_1}{1+\frac{\epsilon_1+\epsilon_2}{2}}\left(\tilde
H-\frac{\epsilon_1-\epsilon_2}{2}K^{(2) }\right)
$$Thus the application of theorem \ref{main} and corollary
\ref{c.main} immediately gives the result.  \qed

\begin{remark}
\label{fre.2}
The ratio between the two frequencies of motion on $\T_\epsilon$ is
$(\epsilon_1-\epsilon_2)/2$. Since such a quantity varies in an
uncountable set, the flow is dense on many of the tori
$\T_\epsilon$. But there are also some tori on which
$\epsilon_1-\epsilon_2$ is rational, thus we have also proved the
persistence of some resonant tori, a phenomenon that does not occur
when applying KAM theory.
\end{remark}

\subsection{A beam vibrating in 2 directions} \label{beam}

As a second example we consider the system of coupled equations
\begin{eqnarray}
u_{tt}+u_{xxxx}-mu=-u f(u^2+v^2)&&\nonumber \\ &&x\in\T \label{bea.1} \\
v_{tt}+v_{xxxx}-mv=-v f(u^2+v^2)&&\nonumber
\end{eqnarray}
This can be thought as a model of a beam which can oscillate in the
$y$ and in the $z$ directions. $u$ represents the displacement in the
$y$ direction and $v$ the displacement in the $z$ directions. In the
case of negative $m$ the terms $mu$ and $mv$ can be interpreted as
centrifugal forces due to the fact that the beam is actually rotating
with a constant velocity around its axis.  To be determined here we
will just consider the case $f(s)=s$.

\eqref{bea.1} is a Hamiltonian system with Hamiltonian function
\begin{equation}
\label{bea.2}
K^{(1)}(U,V,u,v)=\int_{-\pi}^\pi\left[\frac{U^2+u_{xx}^2+mu^2+V^2+v_{xx}^2
    +mv^2}{2} + \frac{\left(u^2+v^2\right)^2}{4}\right] dx\ ,
\end{equation}
where $U$ is the momentum conjugated to $u$ and $V$ the momentum
conjugated to $v$.  There are two symmetries, namely the translations
in the torus and the rotations in the plane $yz$. They are generated
by the two Hamiltonians
\begin{equation}
\label{bea.4}
K^{(2)}(U,V,u,v)=\int_{-\pi}^\pi (U_xu+V_xv)dx\ ,
\quad K^{(3)}(U,u)=\int_{-\pi}^\pi (Uv-Vu)dx\ .
\end{equation}

We have to introduce suitable coordinates. First introduce the Fourier
variables $(U_j,V_j,u_j,v_j)$ as in the case of the nonlinear wave
equation (with $\omega_j=\sqrt{j^4+m}$), then make the coordinate
transformation (for $j\not=0$)
\begin{eqnarray}
\label{bea.8}
u_j=\frac{P_{1,j}+Q_{1,-j}}{\sqrt 2}\ ,\quad
  U_j=\frac{P_{1,-j}-Q_{1,j}}{\sqrt 2} \\
  v_j=\frac{P_{2,j}+Q_{2,-j}}{\sqrt 2}\ ,\quad
  V_j=\frac{P_{2,-j}-Q_{2,j}}{\sqrt 2} \ .
\end{eqnarray}
which transforms the Hamiltonians to
\begin{eqnarray*}
K^{(1)}&=&\sum_{j\in\Z}\omega_j\frac{P_{1,j}^2+Q_{1,j}^2+P_{2,j}^2+
  Q_{2,j}^2}{2} +F(u,v)\ ,
\\
K^{(2)}&=&\sum_{j\in\Z}j\frac{P_{1,j}^2+Q_{1,j}^2+P_{2,j}^2+
  Q_{2,j}^2}{2}
\\
K^{(3)}&=&\sum_{j\in\Z}P_{1,j}Q_{2,j}-P_{2,j}Q_{1,j}\ ,
\end{eqnarray*}
where we denoted
\begin{equation}
\label{bea.8a}
F(u,v)=\int_{-\pi}^\pi \frac{\left(u^2+v^2\right)^2}{4} dx\ .
\end{equation}
Then for all $j$'s make the further change of variables
\begin{eqnarray}
\label{bea.8b}
Q_{1,j}=\frac{p_{1,j}+q_{2,j}}{\sqrt 2}\ ,\quad
  P_{1,j}=\frac{p_{2,j}-q_{1,j}}{\sqrt 2} \\
  Q_{2,j}=\frac{p_{2,j}+q_{1,j}}{\sqrt 2}\ ,\quad
  P_{2,j}=\frac{p_{1,j}-q_{2,j}}{\sqrt 2}\ ,
\end{eqnarray}
which gives the Hamiltonians the form
\begin{eqnarray}
K^{(1)}&=&\sum_{j\in\Z}\omega_j
  \left(\frac{p_{1,j}^2+q_{1,j}^2}2+\frac{p_{2,j}^2+
  q_{2,j}^2}{2} \right) +F(u,v)\ , \\
  K^{(2)}&=&\sum_{j\in\Z}j\left(\frac{p_{1,j}^2+q_{1,j}^2}2+\frac{p_{2,j}^2+
  q_{2,j}^2}{2}\right) \\
  K^{(3)}&=&\sum_{j\in\Z}\frac{p_{1,j}^2+q_{1,j}^2}2-\frac{p_{2,j}^2+
  q_{2,j}^2}{2}
\end{eqnarray}
Then we proceed as in the wave equation: we scale the variables, by
$p=\sqrt\mu \tilde p, q=\sqrt\mu \tilde q$ and we omit tildes.

We fix now three modes that will be non vanishing on the invariant
torus of the linearized system.

Just to give an example we fix the modes $(p_{1,j},q_{1,j})$ with
$j=1,2,3$. We introduce action angle variables $(I_j,\phi_j)$,
$j=1,2,3$ for such modes and look for the linear combinations of the
Hamiltonians having the form \eqref{e.4}. They are given by
\begin{eqnarray}
\label{bea.h1}
H^{(1)}&=&\frac{K^{(1)}+(\omega_2-\omega_3)K^{(2)}+(2\omega_3-3\omega_2)
  K^{(3)}}{\omega_3-2\omega_2+\omega_1} \\
H^{(2)}&=&\frac{-2K^{(1)}+(\omega_3-\omega_1)K^{(2)}+(3\omega_1-\omega_3)
  K^{(3)}}{\omega_3-2\omega_2+\omega_1}
\\
  H^{(3)}&=&\frac{K^{(1)}+(\omega_1-\omega_2)K^{(2)}+(\omega_2-2\omega_1)
  K^{(3)}}{\omega_3-2\omega_2+\omega_1}
\end{eqnarray}
Assumption (A.2) holds with $d=2$.  In order to apply theorem
\ref{main} we choose $\bn=\be=(1,0,0)$. We first verify that the
nonresonance condition (A.3) is fulfilled.
\begin{lemma}
\label{bea.l} {For any $\gamma>0$ small enough} there exists a
subset $\M$ of $(-1,4)$ whose complement has measure $O(\sqrt\gamma)$
such that, if $m\in\M$ then, for any $R>0$, the complementary of set
$\Nc\cap B_R$ has measure estimated {by} $CR\sqrt\gamma$.
\end{lemma}
\proof Instead of \eqref{e.15} we will multiply the quantity $k-\tilde
\Omega_j$ by $\omega_3-2\omega_2+\omega_1$. Thus the quantities that
have to be estimated are
\begin{eqnarray*}
f_{jk}(m,\epsilon)&:=&(1+\epsilon_1-2\epsilon_2+\epsilon_3)\omega_j \\
&&+[(1-\epsilon_2-\epsilon_3)\omega_2
-(1+\epsilon_1+\epsilon_2)\omega_3+(\epsilon_3-\epsilon_2)\omega_1 ] j
\\ &&\pm
[(3\epsilon_2-2\epsilon_3)\omega_1-(1+\epsilon_1-\epsilon_2+\epsilon_3)
\omega_2 +(2+2\epsilon_1-\epsilon_2)\omega_3] \\
&&-k(\omega_3-2\omega_2+\omega_1)
\end{eqnarray*}
where the sign + holds for $\tilde \Omega_{1,j}$, while the sign -
holds for $\tilde \Omega_{2,j}$.

In the particular case
$\epsilon=0$ we thus have to estimate
\begin{equation}
\label{bea.44}
k\omega_1-(2k+j\pm3)\omega_2+(k+j\mp2)\omega_3-\omega_j\ ,
\end{equation}
where the first sign holds for
$\tilde\Omega_{1,j}=\Omega_{1,j}^{(1)}$, while the second holds for
$\tilde\Omega_{2,j}=\Omega_{2,j}^{(1)}$.  Consider \eqref{bea.44} as a
function of $m$. It is easy to see that it vanishes identically only
if (i) $(j,k)=(1,1)$ and the sign is + (which corresponds to $I_1$),
if (ii) $(j,k)=(2,0)$ and the sign is + (which corresponds to $I_2$),
if (iii) $(j,k)=(3,0)$ and the sign is + (which corresponds to
$I_3$). Thus all these cases are excluded. In all the other cases the
function \eqref{bea.44} is nontrivial.  This remains true when
$\epsilon$ is in a small neighborhood of zero.

Notice now that $f_{jk}$ can vanish (or be
small) only if
\begin{equation}
\label{bea.e}
C_1j^2\leq |k|\leq C_2j^2
\end{equation}
Thus we consider
only such a case. First fix a large $J_*$, then exploiting the fact
that $d\omega_j/dm=(2\omega_j)^{-1}$ one has
\begin{equation}
\label{bea.e1}
\left|\frac{\partial f_{jk}}{\partial m}\right|\geq C|k|\geq C'j^2\ .
\end{equation}
It follows that there is at most a segment of length
$C\gamma/|j|^{\tau-2}$ in the $m$ space where $f_{jk}$ is smaller than
$\gamma/|j|^\tau$. Summing over $k$, taking into account \eqref{bea.e}
one gets that except for an interval of masses $m$ of length
$C\gamma/|j|^{\tau}$ one has
$$|f_{jk}|\geq \frac{\gamma}{|j|^{\tau}} \quad \forall k\in\Z
$$
Summing over $j$ one gets that, provided $\tau>1$ one has
$$|f_{jk}|\geq \frac{\gamma}{|j|^{\tau}} \quad \forall k\in\Z, \forall
|j|\geq J_*
$$ except for the $m's$ in a set of measure $C\gamma$. And this is
true {\it for any value of $\epsilon$ small enough}.

It remains to consider the set of $j's$ smaller than $J_*$. But these
are a finite number. Therefore excluding at most a further set of
measure of order $\gamma$ one gets that $\forall \epsilon$ fixed and
close to the origin, $f_{jk}$ fulfills the wonted estimate except for
a set of $m$'s having measure of order $\gamma$. Using Fubini's
theorem one then concludes the proof.\qed

Finally we have to verify the non degeneracy assumption. First
{notice} that $F_{\bn}$ is proportional to $F$,
and that going through the change of coordinates one has
\begin{eqnarray}
\label{bea.301}
u(a,\be t,0,0)=-\sqrt{a_1}\frac{\cos(x+t)}{\sqrt\pi}
-\sqrt{a_2}\frac{{\cos(2x)}}{\sqrt{\pi
}}-\sqrt{a_3}\frac{\cos(3x)}{\sqrt\pi} \\ v(a,\be
t,0,0)=\sqrt{a_1}\frac{\sin(x+t)}{\sqrt{\pi}}
+\sqrt{a_2}\frac{\sin(2x)}{ \sqrt{\pi}} +\sqrt{a_3}\frac{\sin(3x)
}{\sqrt\pi}\ .
\end{eqnarray}
Thus substituting in $F$ and computing the integrals over $x$ and $t$
one obtains that $\langle F_{\bn}\rangle$ is proportional to
$$
a_1^2+a_2^2+a_3^2+4a_1a_2+4a_2a_3+4a_1a_3
$$
which is a non degenerate quadratic form. Therefore also assumption
(A.4) holds.

Applying theorem \ref{main} one thus gets existence of quasiperiodic
solutions with three frequencies.

To get a precise result consider the set
\begin{equation}
\label{bea.16}
\A_*:=\bigcup_{A\in (0,\infty)^3 }\left\{-\frac{\partial
\langle F_{\be}\rangle}{\partial a}(A) \right\}
\end{equation}
and, for $\epsilon\in \A_*$, denote by $A_*(\epsilon)$ the unique
solution of
$$
\epsilon =-\frac{\partial
\langle F_{\be}\rangle}{\partial a}(A) \ ,
$$
and by $\T_{\epsilon,*}$ the torus $p=q=0, I=A_*, \phi\in\T^3$. Thus
the following theorem holds.
\begin{theorem}
\label{bea.teo}
There exists a positive $\mu_*$ with the following property: for any
$\epsilon\in \Nc\cap B_{\mu_*}\cap \A_*$ there exists a unique invariant torus
$\T_{\epsilon}$ such that
\begin{itemize}
\item[1.] the flow on $\T_{\epsilon}$ has frequencies
\begin{eqnarray}
\label{fre.bea}
\left(\frac{\omega_3-2\omega_2+\omega_1}{1+\epsilon_1-2\epsilon_2+\epsilon_3}
,\ - \frac{(1+\epsilon_1)\omega_2-(1+\epsilon_1-
\epsilon_2)\omega_3-\epsilon_2\omega_3}
{1+\epsilon_1-2\epsilon_2+\epsilon_3 },\right.
\\ \nonumber \left.
-\frac{(2+2\epsilon_1-\epsilon_3)\omega_3-
(3+3\epsilon_1-\epsilon_3)\omega_2 +(3\epsilon_2-2\epsilon_3)\omega_1
} {1+\epsilon_1-2\epsilon_2+\epsilon_3 } \right)
\end{eqnarray}
\item[2.]  one has
\begin{equation}
\label{bea.33}
d(\T_{\epsilon},\T_{\epsilon,*})\leq C \norma{A_*(\epsilon)}^2
\end{equation}
\end{itemize}
\end{theorem}
\begin{remark}
\label{mes}
Since the set of $\epsilon\in\Nc\cap B_\mu\cap\A_*$ has positive
measure, one has that for almost all allowed values of $\epsilon$, the
flow is ergodic on $\T_{\epsilon}$.
\end{remark}


\begin{thebibliography}{COR02}

\bibitem[Bal06]{Bal06}
P.~Baldi, \emph{Quasi-periodic solutions of the equation {$v\sb {tt}-v\sb
  {xx}+v\sp 3=f(v)$}}, Discrete Contin. Dyn. Syst. \textbf{15} (2006), no.~3,
  883--903.

\bibitem[Bam00]{bam00}
D.~Bambusi, \emph{Lyapunov center theorem for some nonlinear {PDE}s: a simple
  proof}, Ann. S.N.S. Pisa Cl. Sci \textbf{29} (2000), 823--837.

\bibitem[BG02]{BG02}
D.~Bambusi and G.~Gaeta, \emph{On persistence of invariant tori and a theorem
  by {N}ekhoroshev}, Math. Phys. Electron. J. \textbf{8} (2002), Paper 1, 13
  pp. (electronic).

\bibitem[Bou98]{Bo98}
J.~Bourgain, \emph{Quasi-periodic solutions of {H}amiltonian perturbations of
  2{D} linear {S}chr\"odinger equations}, Ann. of Math. (2) \textbf{148}
  (1998), no.~2, 363--439.

\bibitem[Bou05]{Bou03}
J.~Bourgain, \emph{Green's function estimates for lattice {S}chr\"odinger
  operators and applications}, Annals of Mathematics Studies, vol. 158,
  Princeton University Press, Princeton, NJ, 2005.

\bibitem[BP02]{BP02}
D.~Bambusi and S.~Paleari, \emph{Families of periodic orbits for some {PDE}'s
  in higher dimensions}, Commun. on Pure and Applied Analysis \textbf{1}
  (2002), 269--279.

\bibitem[BP06]{BerPro06}
M.~Berti and M.~Procesi, \emph{Quasi-periodic solutions of completely resonant
  forced wave equations}, Comm. Partial Differential Equations \textbf{31}
  (2006), no.~4-6, 959--985.

\bibitem[BV02]{BV02}
D.~Bambusi and D.~Vella, \emph{Quasi periodic breathers in {H}amiltonian
  lattices with symmetries}, Discrete Contin. Dyn. Syst. Ser. B \textbf{2}
  (2002), no.~3, 389--399.

\bibitem[COR02]{COR02}
P.~Chossat, J.P. Ortega, and T.S. Ratiu, \emph{Hamiltonian {H}opf bifurcation
  with symmetry}, Arch. Ration. Mech. Anal. \textbf{163} (2002), no.~1, 1--33.

\bibitem[CW93]{CW92}
W.~Craig and C.~E. Wayne, \emph{Newton's method and periodic solutions of
  nonlinear wave equations}, Comm. Pure Appl. Math. \textbf{46} (1993),
  1409--1498.

\bibitem[CY00]{CY00}
L.~Chierchia and J.~You, \emph{{KAM} tori for 1{D} nonlinear wave equations
  with periodic boundary conditions}, Comm. Math. Physics \textbf{211} (2000),
  497--525.

\bibitem[EK06]{EK06}
H.~L. Eliasson and S.~B. Kuksin, \emph{{KAM} for non-linear {S}chroedinger
  equation}, Annals of Math., to appear. Preprint 2006.

\bibitem[Kuk87]{K87}
S.~B. Kuksin, \emph{Hamiltonian perturbations of infinite-dimensional linear
  systems with an imaginary spectrum}, Funct. Anal. Appl. \textbf{21} (1987),
  192--205.

\bibitem[Kuk92]{Kuk92}
S.B. Kuksin, \emph{An infinitesimal {L}iouville-{A}rnold theorem as a criterion
  of reducibility for variational {H}amiltonian equations}, Chaos Solitons
  Fractals \textbf{2} (1992), no.~3, 259--269.

\bibitem[Kuk00]{K2}
S.~B. Kuksin, \emph{Analysis of {H}amiltonian {PDEs}}, Oxford University Press,
  Oxford, 2000.

\bibitem[Mon97]{Mon97}
J.~Montaldi, \emph{Persistence and stability of relative equilibria},
  Nonlinearity \textbf{10} (1997), no.~2, 449--466.

\bibitem[Nek94]{N94}
N.~N. Nekhoroshev, \emph{The {P}oincar\'e-{L}yapunov-{L}iouville-{A}rnold
  theorem}, Funktsional. Anal. i Prilozhen. \textbf{28} (1994), no.~2, 67--69.

\bibitem[Pro05]{Pro05}
M.~Procesi, \emph{Quasi-periodic solutions for completely resonant non-linear
  wave equations in 1{D} and 2{D}}, Discrete Contin. Dyn. Syst. \textbf{13}
  (2005), no.~3, 541--552.

\end{thebibliography}


\providecommand{\bysame}{\leavevmode\hbox to3em{\hrulefill}\thinspace}
\providecommand{\MR}{\relax\ifhmode\unskip\space\fi MR }
\providecommand{\MRhref}[2]{%
  \href{http://www.ams.org/mathscinet-getitem?mr=#1}{#2}
}
\providecommand{\href}[2]{#2}

\end{document}